\documentclass[english,12pt]{article}
\usepackage[T1]{fontenc}
\usepackage{float}
\usepackage{fullpage}
\usepackage{amsmath,amsthm}
\usepackage{amssymb}
\usepackage{graphicx}
\usepackage{subfigure}
\usepackage{color}
\usepackage{epsfig}
\usepackage{epstopdf}
\usepackage{hyperref}
\usepackage{authblk}

\allowdisplaybreaks
\numberwithin{equation}{section}
\title{Should exponential integrators be used for advection-dominated problems?}
\author{Lukas Einkemmer\thanks{\href{mailto:lukas.einkemmer@uibk.ac.at}{lukas.einkemmer@uibk.ac.at}}}
\author{Trung-Hau Hoang\thanks{\href{mailto:trunghaugg@gmail.com}{trunghaugg@gmail.com}}}
\author{Alexander Ostermann\thanks{\href{mailto:alexander.ostermann@uibk.ac.at}{alexander.ostermann@uibk.ac.at}}}
\affil{Department of Mathematics, University of Innsbruck, Austria}
\date{}

\begin{document}

\maketitle

\begin{abstract}\noindent
In this paper, we consider the application of exponential integrators to problems that are advection dominated, either on the entire or on a subset of the domain. In this context, we compare Leja and Krylov based methods to compute the action of exponential and related matrix functions. We set up a performance model by counting the different operations needed to implement the considered algorithms. This model assumes that the evaluation of the right-hand side is memory bound and allows us to evaluate performance in a hardware independent way. We find that exponential integrators perform comparably to explicit Runge--Kutta schemes for problems that are advection dominated in the entire domain. Moreover, they are able to outperform explicit methods in situations where small parts of the domain are diffusion dominated. We generally observe that Leja based methods outperform Krylov iterations in the problems considered. This is in particular true if computing inner products is expensive.
\end{abstract}

\section{Introduction}

In the last few decades, exponential integrators have gained popularity as an option for solving complex systems of differential equations that exhibit stiffness (see, e.g., the review \cite{HO2010}). The main idea behind exponential integrators is as follows: the right-hand side of the problem is linearized, the linear component is integrated exactly and the remaining nonlinearity is discretized in an explicit way. In order to implement such schemes, approximations of the actions of the matrix exponential and related matrix functions are required. Doing this efficiently is important to obtain a competitive method.

In situations where the diffusion/advection coefficient is constant, the domain is a simple rectangle, and either periodic or homogeneous Dirichlet/Neumann boundary conditions are applied, the fast Fourier transform (FFT) or related techniques can be employed (see, e.g., \cite{Karle2006,Crouseilles2018,Crouseilles2020,Caliari2022}). If these assumptions do not hold, alternative methods must be used. Conventional approaches such as Pad\'{e} approximations or diagonalization, however, are only practical when dealing with small systems. This paper focuses on two widely used classes of methods for approximating the action of a matrix functions on vectors in the context of large-scale systems. Specifically, we examine Krylov subspace methods (see, e.g., \cite{e6fd4a4b-1a44-3d13-b09d-a560ff78d9f9,10.1145/2168773.2168781}) and Leja interpolation (see, e.g., \cite{10.1007/s10543-013-0446-0,CALIARI200479}).

Exponential methods have been investigated extensively for diffusion-dominated problems (see, e.g., \cite{HO2010,doi:10.1137/S1064827595295337,Loffeld2013}). In this situation, it can be shown that they exhibit significantly improved performance compared to explicit methods and similar or in some cases even improved performance compared to implicit methods. Most commonly this comparison is performed under the assumption that no or no good preconditioner is available, although some ways to incorporate the solution of a simplified but similar problem into exponential integrators have been explored \cite{Castillo1998,Eshof2006,Einkemmer2023}. Significantly fewer research has been done in the context of advection-dominated problems (see \cite{CALIARI2009568,pub.1146932560,EINKEMMER2017550,doi:10.1137/080717717}), the situation we consider in this paper. Here, we consider linear advection-diffusion problems with non-constant diffusion coefficients and a 2D compressible isothermal Navier--Stokes problem.

To perform the comparisons, we evaluate the costs associated with the different methods. Our goal it to compare performance in a hardware independent way. To accomplish this, we develop a performance model (as is commonly done in computer science \cite{treibig2009,stengel2015,lagraviere2019}) that counts operations such as matrix-vector multiplications, inner products, scalar multiplications, and performing linear combinations. This also allows us, for example, to investigate what effect expensive inner products (as is common in the case of modern distributed memory supercomputers) have on the overall cost.

The paper is structured as follows. In Section \ref{testproblem}, we introduce the test problems for the numerical investigations. The performance model used to evaluate computational cost is explained in Section \ref{computationalmodel}. In Sections \ref{num1d} and \ref{num2d}, we present the numerical experiments for these problems. Finally, Section \ref{sectionconcludesionpaper1} provides the conclusions of the study.

\subsection{Exponential Rosenbrock methods}

In this paper, we compare the performance of explicit Runge--Kutta methods with exponential Rosenbrock schemes. While the former are well-established (see, e.g., the textbooks \cite{hairer2008solving,book}), we recall two representations of the latter that will be used in this investigation. To present the schemes, we consider the autonomous problem
\begin{equation}\label{eq9}
	u^{\prime}(t)=F(u(t)), \quad u(t_0)=u_0,
\end{equation}
where $F$ is a nonlinear function of $u$. Let $u_k$ be the numerical solution, which approximates the exact solution $u(t_k)$ at time $t_k$. By linearising equation \eqref{eq9} in each step at $u_k$, we obtain
\begin{equation}\label{eq11}
	u^{\prime}(t)=J_k u(t) + g_k( u(t)),	
\end{equation}
where
$$
\begin{aligned}
	J_k  =\frac{\partial F}{\partial u}(u_k), \qquad g_k( u(t)) = F(u(t))-J_k u(t).
\end{aligned}
$$
The simplest exponential Rosenbrock method, the so-called exponential Rosenbrock--Euler method (see \cite{HO2010}), has the form
\begin{equation}\label{exponentialRosenbrockEuler}
\begin{aligned}
	u_{k+1} &=u_{k}+\tau \varphi_{1}(\tau {J}_{k}) F(u_{k}),\qquad \varphi_1(z) = \frac{\mathrm{e}^z-1}z,
\end{aligned}	
\end{equation}
where $\tau$ denotes the time step size. This method requires the action of the matrix function $\varphi_1$ on the right-hand side $F$ of the problem, and it is second-order convergent (see \cite{HO2010,doi:10.1137/080717717}). Higher-order methods have been developed in the literature (see, e.g., \cite{doi:10.1137/080717717}). In view of recent developments, two-stage fourth-order exponential integrators for time-dependent partial differential equations (PDEs) have been constructed (for further details, we refer the reader to \cite{LUAN201791}). In the course of our experiments, we employed the so-called $\mathtt{exprb42}$ scheme

\begin{equation}\label{eq1}
\begin{aligned}
	U_{k 2} &= u_{k} + \tfrac{3}{4} \tau  \varphi_{1}\bigl( \tfrac{3}{4} \tau J_k \bigr)F(u_{k}), \\[2pt]
	u_{k+1 }&= u_k +\tau \varphi_1(\tau J_k)F(u_{k}) + \tfrac{32}{9}\tau  \varphi_3(\tau J_k ) \left(g( U_{k2}) - g( u_k)\right),
\end{aligned}	
\end{equation}
where
$$
\varphi_3(z) = \frac{\mathrm{e}^z-1-z -z^2/2}{z^3}.
$$
This scheme is interesting from a computational point of view because it requires only one internal stage at each time step. Additionally, updating  $u_{k+1}$ involves only three actions of matrix functions, a notable reduction compared to the four products in $\mathtt{exprb43}$ (see \cite{HO2010}).

\section{Test problems}\label{testproblem}

This section presents the test problems employed for our investigation. The first problem under consideration is a linear advection-diffusion equation, in which the relative strength of diffusion with respect to advection is controlled through the diffusion coefficient. For this problem, a one-dimensional version is sufficient for the purposes of this study. The second problem are the nonlinear compressible isothermal Navier--Stokes equations, which are a typical representative of an advection-dominated situation.

\subsection{Linear advection-diffusion problem}

\begin{figure}[h!]
\subfigure{\includegraphics[width=0.51\textwidth]{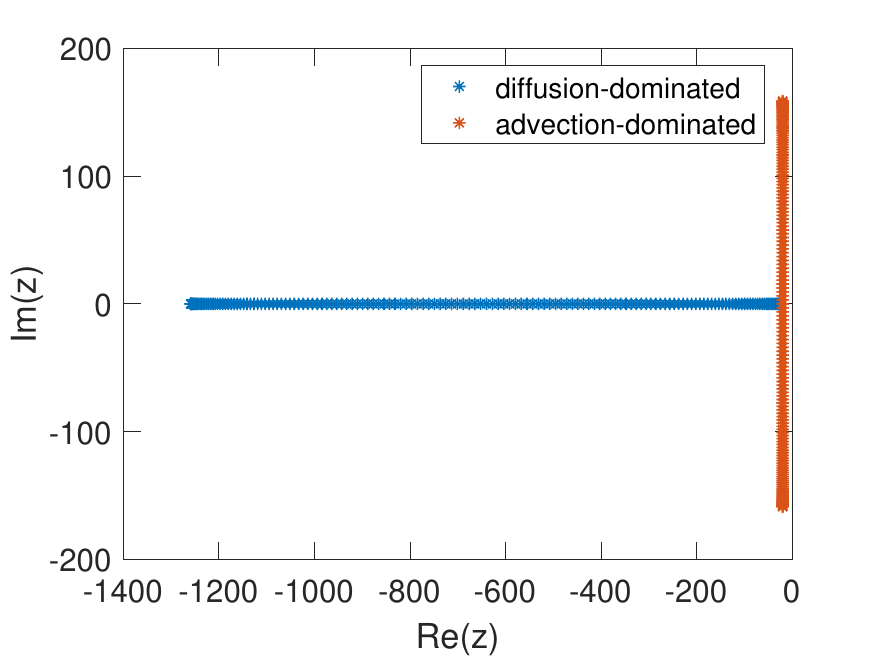}}			
\subfigure{\includegraphics[width=0.51\textwidth]{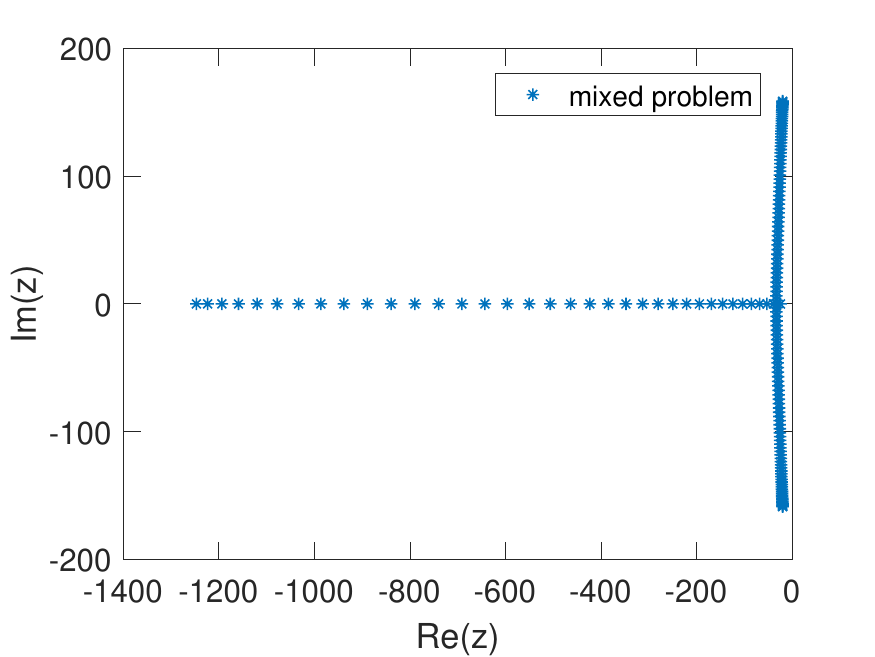}}\label{fig4c}			
\caption{Spectrum of the matrix $-(A_h + B_h) $ for the linear problem \eqref{problem1} for different choices of $\kappa$. In the left figure, the choice $\kappa = 1/80$ indicates a diffusion-dominated problem, while the choice $\kappa = 1/2560$ indicates an advection-dominated problem. The figure on the right, where the function $\kappa$ is defined as described in Section \ref{mixproblem}, shows the spectrum of a mixed problem which, depending on the position, is either diffusion- or advection-dominated.}
\label{spectrum1}
\end{figure}

We consider the following linear advection-diffusion problem
\begin{equation}\label{problem1}
	\begin{aligned}
		& \partial_t  u(t,x) - \kappa(x) \partial_{xx} u (t,x) + \partial_x u(t,x)= 0,\qquad (t,x)\in[0,1]\times[0,1], \\
		& u(0,x) = x(1-x),
	\end{aligned}	
\end{equation}
with homogeneous Dirichlet boundary conditions. The diffusion coefficient $\kappa$ acts as a parameter and is a given (positive) function. For the spatial discretization of~\eqref{problem1}, we consider an equidistant grid with $n$ interior grid points and mesh width $h = 1/(n+1)$. We use standard second-order finite differences for the discretization of $ -\kappa(x) \partial_{xx}$ and $ \partial_{x}$, resulting in the matrices $A_h$ and $B_h$, respectively. We note that a value of $\kappa$ which is larger than $8h$ indicates a diffusion-dominated problem, while a sufficiently low value implies an advection-dominated problem. The influence of $\kappa$ on the spectrum of the matrix $-(A_h+B_h)$ is illustrated in Figure~\ref{spectrum1}.

\subsection{Compressible isothermal Navier--Stokes problem}\label{discussNavierstokes}

\begin{figure}[h!]
\subfigure{\includegraphics[width=0.51\textwidth]{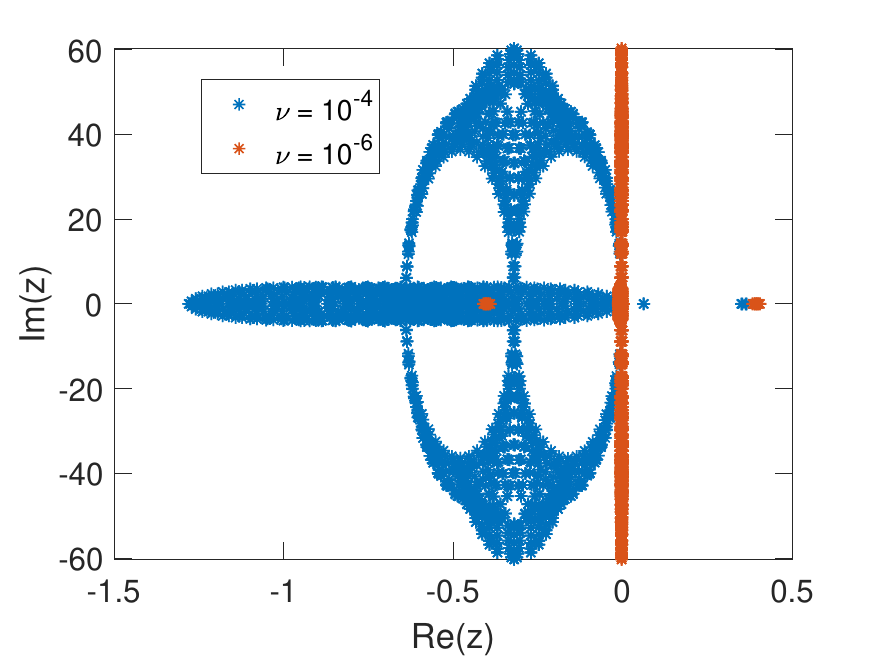}}
\subfigure{\includegraphics[width=0.51\textwidth]{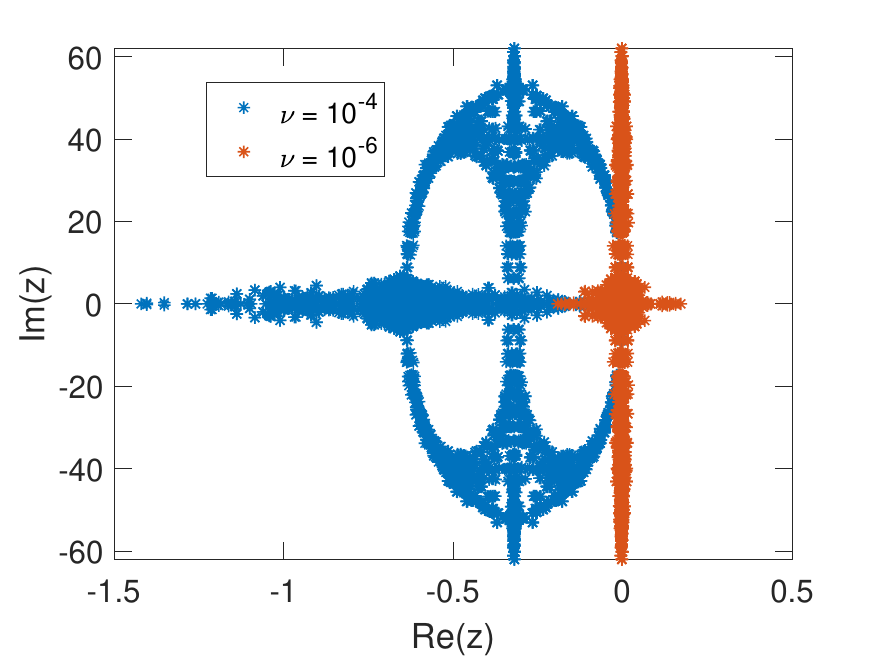}}	
\caption{The spectra of the Jacobians evaluated at the initial value (left figure) and at the final value (right figure) of a numerical solution of the two-dimensional compressible isothermal Navier--Stokes problem \eqref{eq41} with diffusion coefficients of \(\nu = 10^{-4}\) and \(\nu = 10^{-6}\), respectively, indicating advection-dominated scenarios.}
\label{spectrum3}
\end{figure}

In space dimension two, the isothermal Navier--Stokes equations can be written as follows
\begin{equation}\label{eq41}
\begin{aligned}
	&\partial_t \rho + \partial_x (\rho u) + \partial_y (\rho v)  = 0,	\\
	&\partial_{t} u + u \partial_{x} u+v \partial_{y} u + \tfrac{1}{\rho}\partial_{x} \rho =  \nu\left(\partial_{x x} u+\partial_{y y} u\right),	\\
	&\partial_{t} v+u \partial_{x} v+v \partial_{y} v +\tfrac{1}{\rho} \partial_{y} \rho =\nu\left(\partial_{x x} v+\partial_{y y} v\right).
\end{aligned}
\end{equation}
The first equation represents the mass continuity equation, while the other two equations are the Cauchy momentum equations. Here, the diffusion parameter $\nu$ is taken to be constant, and the problem is considered on the domain $\Omega = [0,1]^2$ and for $t\in[0,12]$, subject to periodic boundary conditions. Numerical experiments are performed with initial data corresponding to a shear flow. The above system can also be expressed in the compact form
\begin{equation}
	\label{1729}
	U^{\prime}(t)=F(U(t)), \quad U(0)=U_{0} ,
\end{equation}
where
\begin{equation}\label{eq42}
	U = \left[\begin{array}{l}
		\rho \\
		u \\
		v
	\end{array}\right], \quad F(U) = \left[\begin{array}{c}
	F_1(U) \\
	F_2(U) \\
	F_3(U)
\end{array}\right] =
	\left[\begin{array}{c}
		-\partial_x (\rho u) - \partial_y (\rho v) \\
		- u \partial_{x} u - v \partial_{y} u - \frac{1}{\rho} \partial_{x} \rho + \nu (\partial_{xx}u + \partial_{yy}u) \\
		-u \partial_{x} v-v \partial_{y} v - \frac{1}{\rho} \partial_{y} \rho +\nu\left(\partial_{x x} v+\partial_{y y} v\right)
	\end{array}\right].
\end{equation}
The exponential integrators under consideration are based on a continuous linearization along the numerical solution and thus require the Jacobian. For further technical details on the numerical implementation of the Jacobian matrix and the right-hand side vector $F(U)$, we refer the reader to the Appendix \ref{Appendix1}. The parameter $\nu$ plays a crucial role in the experiments, as it can be adjusted to fix the magnitude of diffusion in the problem. For illustration, we present the spectrum of the Jacobian matrix with $40$ degrees of freedom in each direction and two distinct values of $\nu$ in Figure \ref{spectrum3}.

\section{Performance model}\label{computationalmodel}

This section presents a model that is used for a comparison of the considered numerical methods. The operations commonly employed for the solution of large-scale differential equations include matrix-vector multiplications with sparse or even stencil matrices, inner products, and forming linear combinations of vectors. While these operations are relatively simple, the vectors will usually exceed the capacity of the cache memory due to high dimensionality. This means that the problem is memory-bound. As a consequence, we will have frequent memory accesses and data transfers between the main memory and the cache, which is an expensive task on modern computers.

In order to gain insight into the computational cost, our experiments count the total number of memory operations. Although, ultimately we are interested in the wall clock time, this requires a well optimized code specifically tailored to the hardware architecture under consideration. The advantage of the performance model is that it gives us an indication of the cost of the algorithm that is hardware-independent (see also the discussion in Section \ref{Desktopvssuper}).

\subsection{Performance model for the advection-diffusion problem}\label{subsec:ad}

Given that our problem is memory-bound, the computational cost on a modern computer is determined by the memory operations that are required by the different operations. Let $u$ be a state vector of length $n$, $\alpha u$ its product with a scalar, L$k$ a linear combination of $k$ state vectors, $\left<u ,v \right>$ the inner product of two state vectors, and finally $Mu$ the product of a state vector with a stencil matrix (like the aforementioned matrices $A_h$ or $B_h$). The memory operations involved in computing (and storing these objects) are given in Table \ref{table27}. The total computational cost is simply the sum of the memory operations required for all operations.

\begin{table}[H]
\centering
\begin{tabular}{|c|c|c|c|c|c|c|}
		\hline
		& fetch $u$ & store $u$ & $\alpha u$ & L$k$ & $ \left<u ,v \right>$ & $ Mu $  \\
		\hline
		Memory operations & $n$ & $n$ & $2n$ & $(k+1)n$ & $2n$ & $2n$   \\
		\hline
\end{tabular}
\caption{Cost of memory operations  for the advection-diffusion problem.\label{table27}}
\end{table}

\subsection{Performance model for the Navier--Stokes problem}

In the context of the two-dimensional compressible isothermal Navier--Stokes problem, the vectors $\rho$, $u$, and $v$ have the length $N$. Let $U = [\rho, u, v]^T$, $V$, and $W$ denote state vectors of \eqref{eq42}. These vectors are of length $3N$. In addition to Section~\ref{subsec:ad}, we consider the matrix-vector multiplication $J(U)W$ and the evaluation $F(U)$, both of which are described in Appendix \ref{Appendix1}. In order to compute $J_{i1}(U)w_1+ J_{i2}(U)w_2 + J_{i3}(U)w_3$ for $i=1,2,3$, one needs to fetch the vectors $w_1$, $w_2$, $w_3$, $u$, $v$, $\rho$ and store the result. This requires $7N$ memory operations. In the same way, the computation of $F_i(U)$ requires $4N$ memory operations. This explains the number of memory operations given in Table \ref{table28}.

\begin{table}[H]
	\centering
	\begin{tabular}{|c|c|c|c|c|c|c|c|c|c|}
		\hline
		& fetch $U$ & store $U$ & $\alpha U$ & L$k$ &  $ \left<U ,V \right>$ & $ J(U)W $ & $F(U)$ \\
		\hline
		Memory operations & $3N$ & $3N$ & $6N$ & $3(k+1)N$ & $6N$ & $21N$ & $12N$ \\
		\hline
	\end{tabular}
	\caption{Cost of memory operations for the two-dimensional Navier--Stokes problem.\label{table28}}
\end{table}%

\subsection{ Desktop computers and supercomputers  }\label{Desktopvssuper}

The number of memory operations necessary for the execution of an inner product varies significantly across different machines, particularly supercomputers, due to the increased communication demands associated with such operations. Consequently, this operation is inexpensive on a desktop computer but expensive on a supercomputer. To account for this discrepancy, we propose to multiply the number of inner products by a factor $\zeta$, which is determined by the network speed, in order to differentiate between desktop computers and supercomputers. For the purposes of our investigation, we set $\zeta = 1$ for desktop computers and $\zeta = 10$ as a representative value for supercomputers.

\section{Numerical results for linear diffusion- and advection-dominated problems}\label{num1d}

In this section, we undertake a comparative analysis of the computational cost associated with solving the test problems either by exponential integrators or the classical explicit Runge--Kutta methods RK2 and RK4. In the case of a linear autonomous problem, the only requisite for exponential integrators is the computation of the action of a matrix exponential. For this task, we are employing either a Krylov subspace method or Leja interpolation. For Krylov we compute the orthogonal basis by the Arnoldi process and then compute the $\varphi_i$ function for the obtained Hessenberg matrix (see, e.g., \cite{HO2010}. We employ the standard termination criterium in \cite{10.1145/2168773.2168781}. For Leja we use a precomputed set of Leja points in the interval $[-2,2]$ which are then scaled and shifted to the spectrum of the matrix (see, e.g., \cite{CALIARI2009568}). As termination criterium the $L^2$ norm of the last computed term is used.
The aim of this study is to compare the performance of the aforementioned methods in scenarios where diffusion and/or advection play a dominant role. We recall that for values of the parameter $\kappa$ which are larger than $8h$, the problem is diffusion-dominated problem, while for small values it is advection dominated. In the numerical experiments conducted on problem \eqref{problem1}, we use $n=159$ interior grid points.

\begin{figure}[h!]
\subfigure{\includegraphics[width=0.51\textwidth]{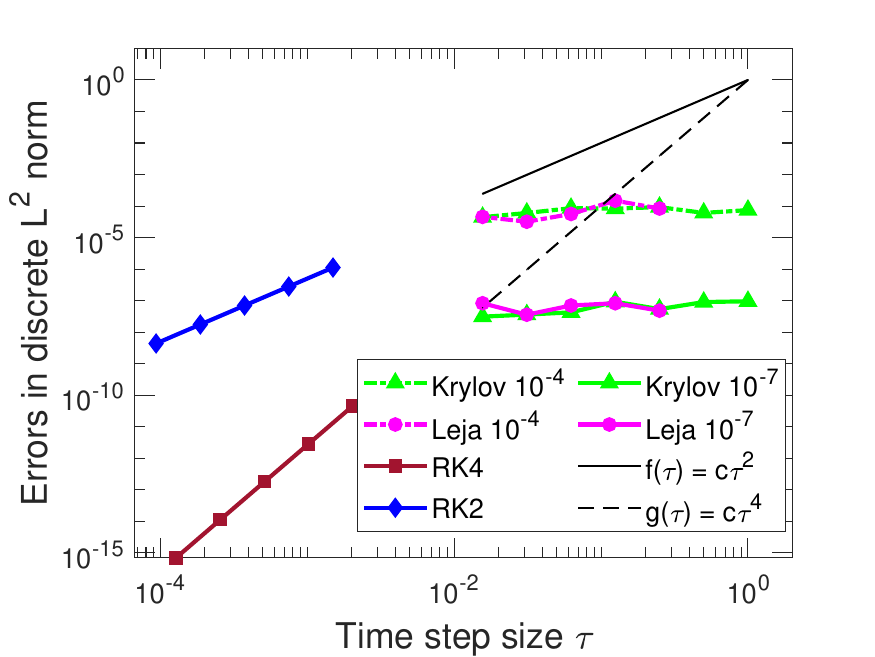}\label{fig1a}}
\subfigure{\includegraphics[width=0.51\textwidth]{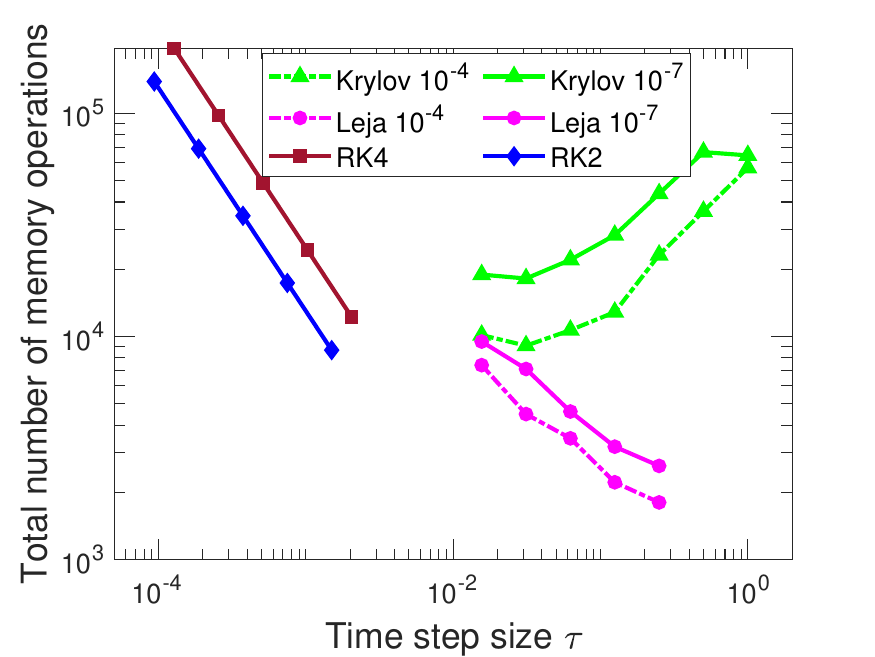}\label{fig1b}}	
\subfigure[ $\tau = \frac{1}{4}$ ]{\includegraphics[width=0.51\textwidth]{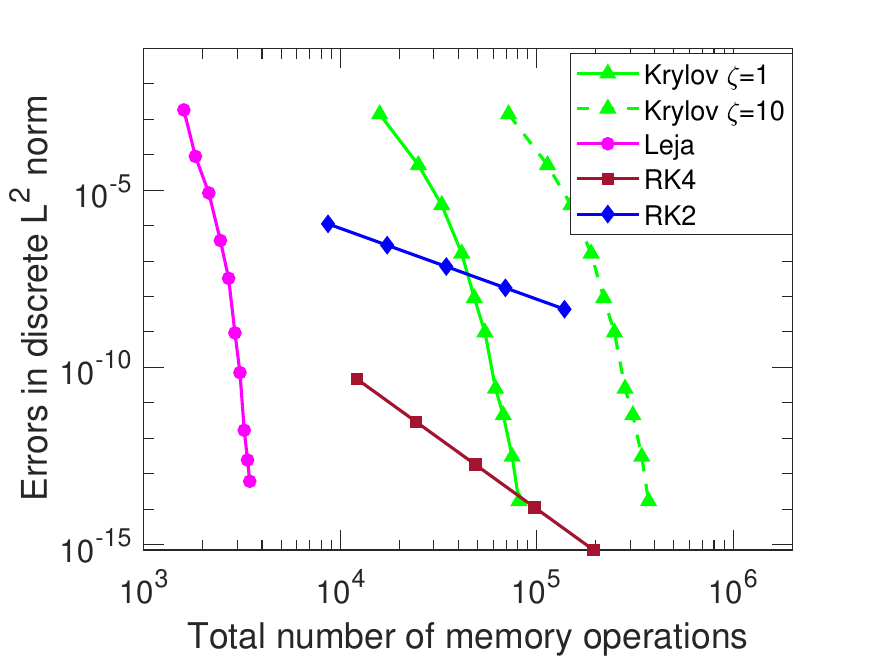}\label{fig1c}}
\subfigure[ $\tau = \frac{1}{64}$]{\includegraphics[width=0.51\textwidth]{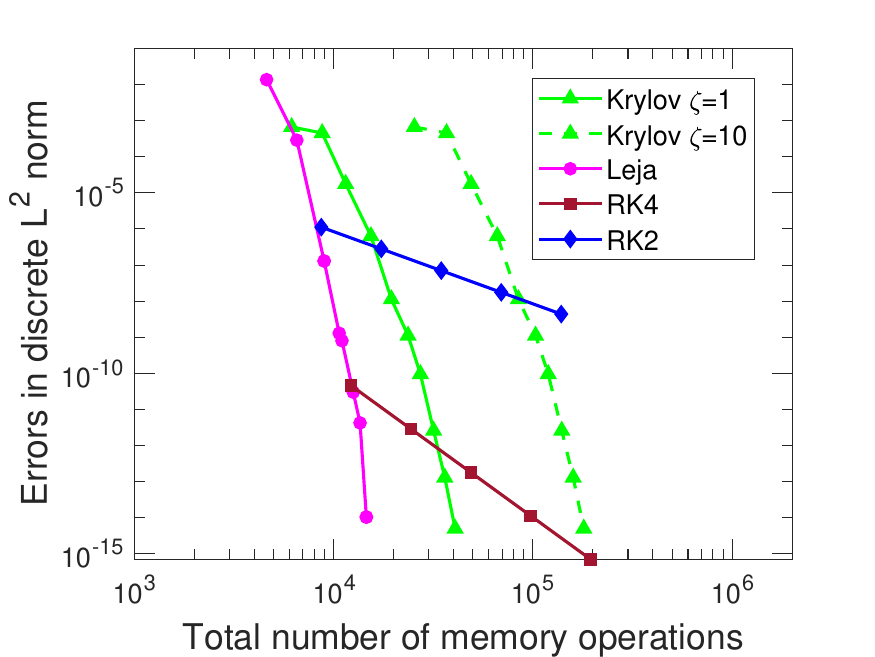}\label{fig1d}}		
\caption{The numerical results for the diffusion-dominated problem \eqref{problem1} with $\kappa=\frac1{80}$ are presented in the upper two panes. The figures illustrate the achieved accuracy and computational cost of the considered methods for evaluating $\exp(- (A_h+B_h))u_0$ as a function of the time step size. The tolerance is chosen such that the achieved accuracy for the exponential integrators is $10^{-4}$ (dashed-dotted line) and $10^{-7}$ (solid line), respectively. The lower two panes depict the achieved accuracy of the considered methods as a function of the computational cost for two chosen values of $\tau$ and $\zeta$. \label{fig1}}
\end{figure}

\subsection{Diffusion-dominated problem}\label{diffusiondominated}

We set $\kappa = \frac{1}{80}$ to ensure that the problem falls within the diffusion-dominated regime. In order to employ the Krylov iteration or Leja interpolation, it is necessary to prescribe both a time step size and a desired accuracy. In the course of our experiments, the prescribed tolerances were $10^{-4}$ and $10^{-7}$. The results of the experiments are collected in Fig.~\ref{fig1}. From this figure it is clear that the Krylov and Leja methods offer several advantages. These methods permit significantly larger step sizes, up to 500 times larger for Krylov and 125 times larger for Leja, compared to the step sizes required by RK2 and RK4. For large time steps, an additional substepping procedure might be necessary, which does not provide any performance advantage. Therefore, we have chosen $\tau = \frac{1}{4}$ as the maximum time step size for Leja.

The Leja method shows superior performance in terms of computational cost in comparison to explicit Runge--Kutta methods. Furthermore, it is also more efficient than Krylov subspace methods, as it necessitates only matrix-vector products, which are less computationally demanding than the operations required by the former. Krylov subspace methods, which employ the Arnoldi process to construct an orthonormal basis, show less perfomance, particularly at large step sizes. This is due to the fact that the use of large time step requires a long recurrence for constructing the orthonormal basis and the upper Hessenberg matrix. This disadvantage is mitigated by using smaller step sizes (resulting in a shorter recurrence). On the other hand, Krylov methods require the evaluation of numerous inner products, which represents a significant disadvantage on modern supercomputer installations in comparison to the other methods under consideration. The plot for $\zeta=10$ in our figures clearly shows this disadvantage of Krylov methods.

\begin{figure}[t!]
\subfigure{\includegraphics[width=0.51\textwidth]{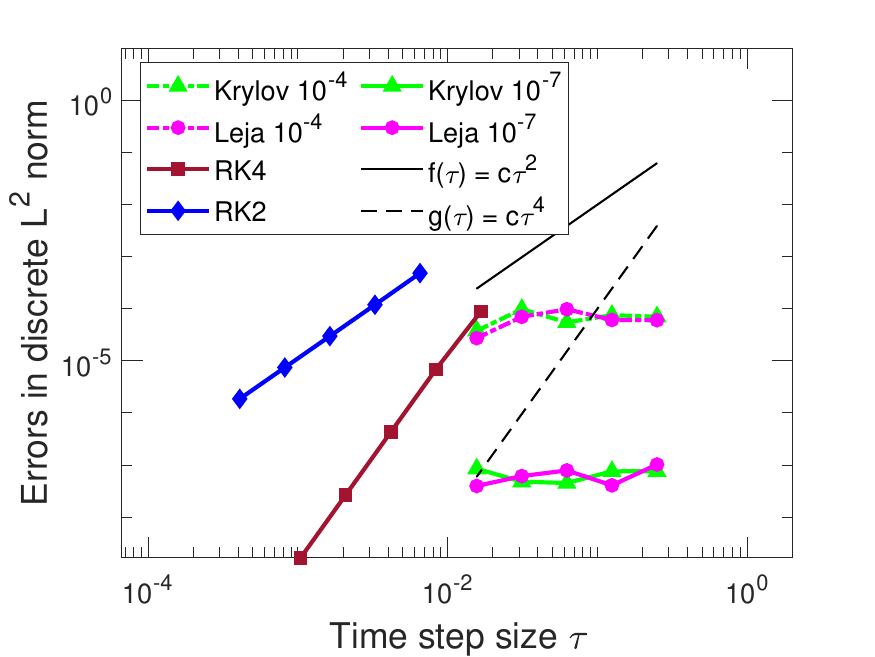} \label{fig2a}}
\subfigure{\includegraphics[width=0.51\textwidth]{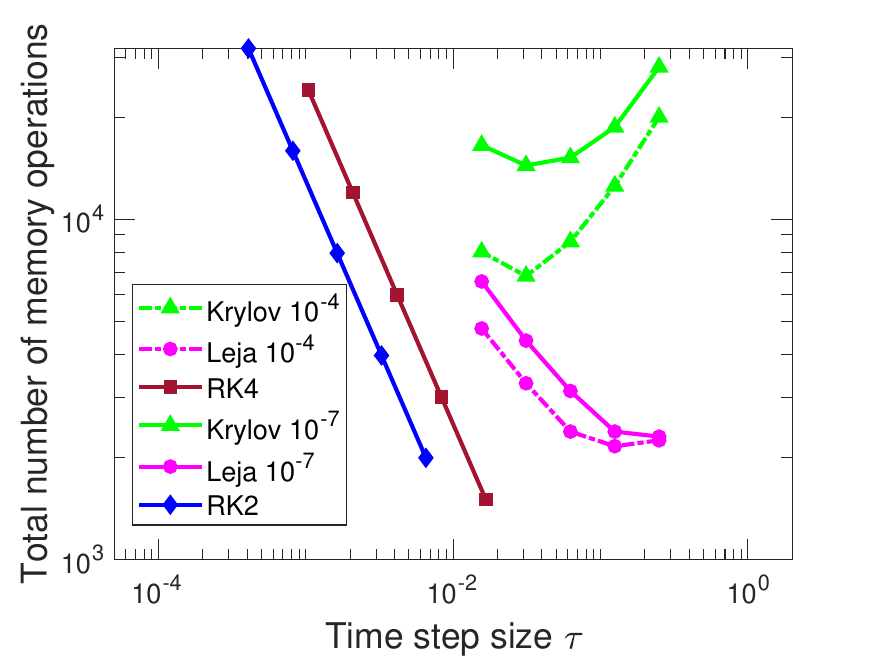} \label{fig2b}}	
\subfigure[$\tau = \frac{1}{4}$ ]{\includegraphics[width=0.51\textwidth]{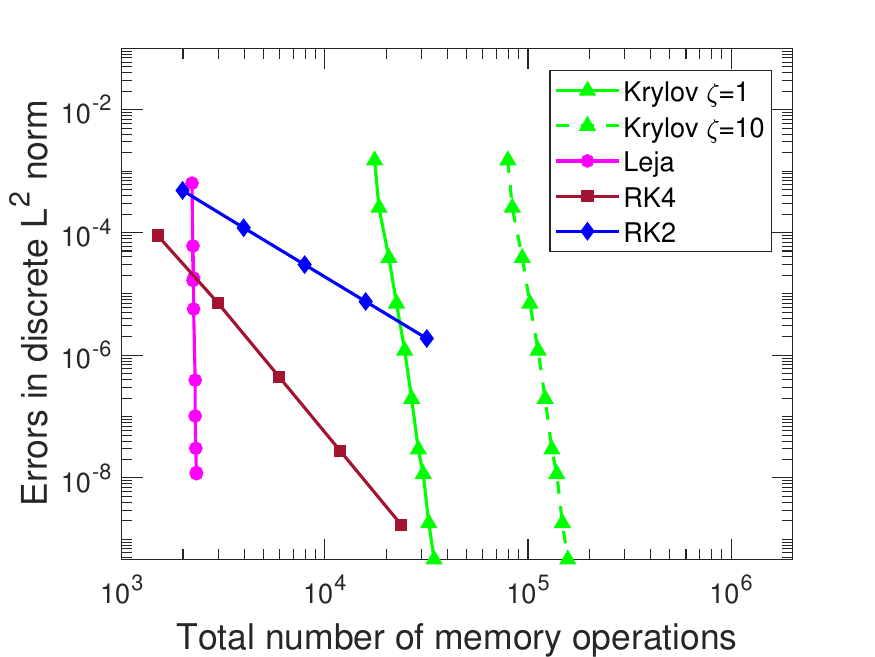} \label{fig2c}}
\subfigure[ $\tau = \frac{1}{64}$ ]{\includegraphics[width=0.51\textwidth]{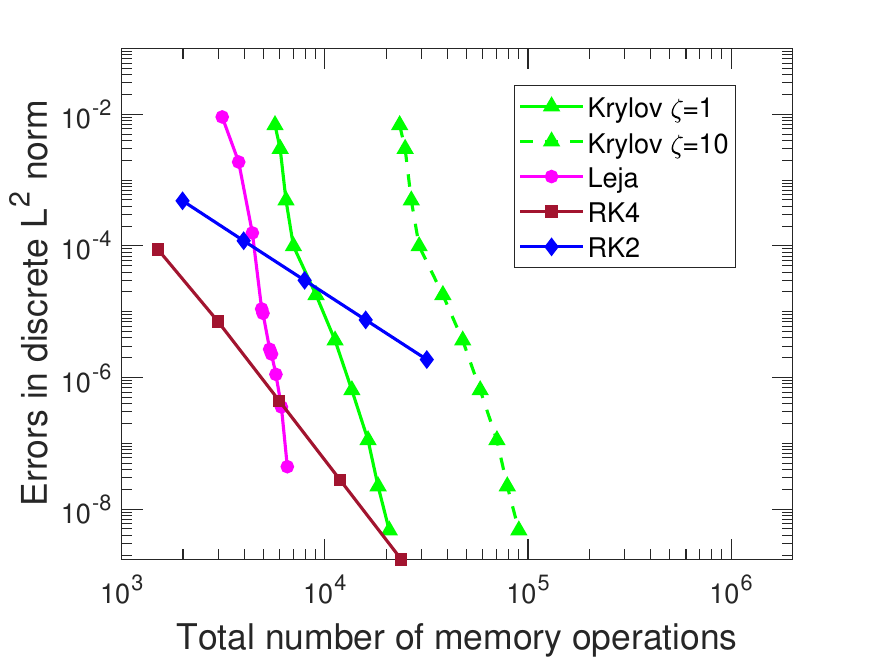} \label{fig2d}}					
\caption{The numerical results for the advection-dominated problem \eqref{problem1} with $\kappa=\frac1{2560}$ are shown in the upper two panes. The figures illustrate the achieved accuracy and computational cost of the considered methods for evaluating $\exp(- (A_h+B_h))u_0$ as a function of the time step size. The tolerance is chosen such that the achieved accuracy for the exponential integrators is $10^{-4}$ (dashed-dotted line) and $10^{-7}$ (solid line), respectively. The lower two panes depict the achieved accuracy of the considered methods as a function of the computational cost for two different values of $\tau$ and $\zeta$. \label{fig2}}
\end{figure}

For a medium accuracy of $10^{-4}$, the action of the exponential function can be approximated by the Leja method with a reasonable computational cost, whereas explicit methods are unable to reach this regime. For a tolerance of $10^{-7}$, the Leja method also outperforms traditional integrators such as RK2 and RK4. As we can see from Figure \ref{fig1}, Krylov and Leja methods demonstrate the capability to achieve much higher accuracy without a significant increase in computational effort, distinguishing them from explicit Runge--Kutta methods. In conclusion, the Leja method is able to achieve any desired level of accuracy while requiring fewer computational resources than RK2 and RK4 methods.

\subsection{Advection-dominated problem}\label{advectiondominated}

Following the promising experiments with the diffusion-dominated problem, we continue to investigate an advection-dominated scenario. To ensure that the problem is in the advection-dominated regime, we choose $\kappa = \frac1{2560}$. For both the Krylov and Leja methods, we select the maximum time step size to be equal to $\frac14$.

The results of the experiments are presented in Figure \ref{fig2}. Exponential methods permit larger step sizes than RK2 and RK4, up to a factor of 38 and 15, respectively. In addition, the Leja method performs competitively in terms of computational cost compared to explicit approaches. However, we observe that when substeps are used in the Leja method, they are more expensive, as in the diffusion-dominated case. On the other hand, the Krylov method benefits from substepping in time. Furthermore, the exponential methods exhibit the capability to achieve significantly higher accuracy without a substantial increase in computational cost, distinguishing them from explicit Runge--Kutta methods.

\subsection{Mixed problem}\label{mixproblem}

\begin{figure}[t!]
\subfigure[]{\includegraphics[width=0.51\textwidth]{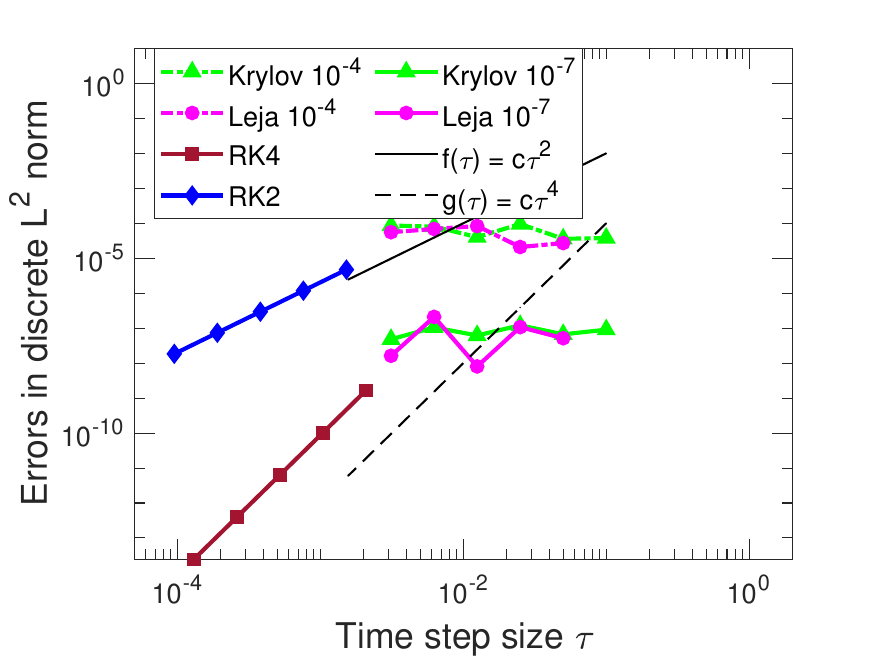} \label{fig4a}}
\subfigure[]{\includegraphics[width=0.51\textwidth]{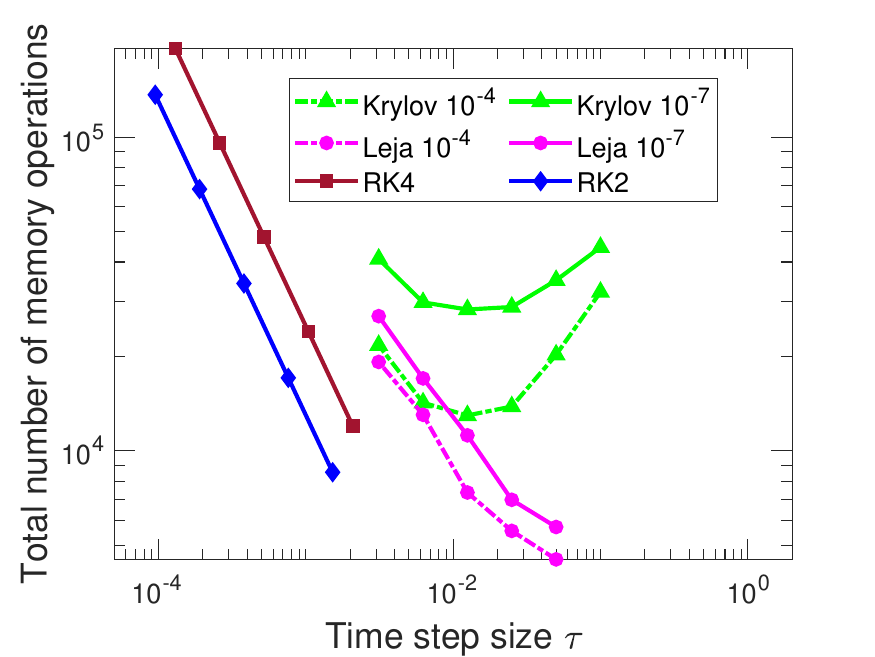} \label{fig4b}}	
\subfigure[] {\includegraphics[width=0.51\textwidth]{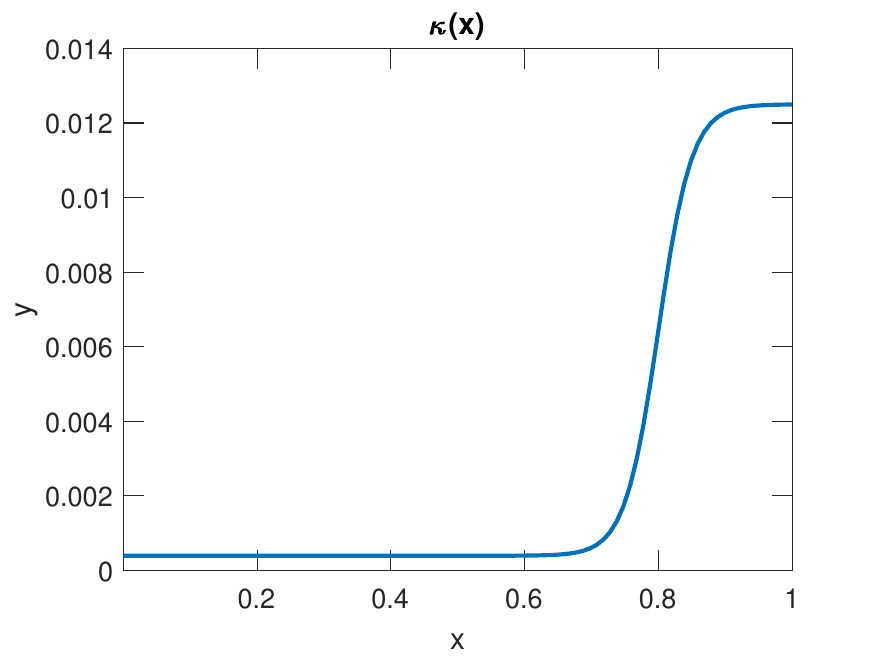} \label{fig4cc}}
\subfigure[  ]{\includegraphics[width=0.51\textwidth]{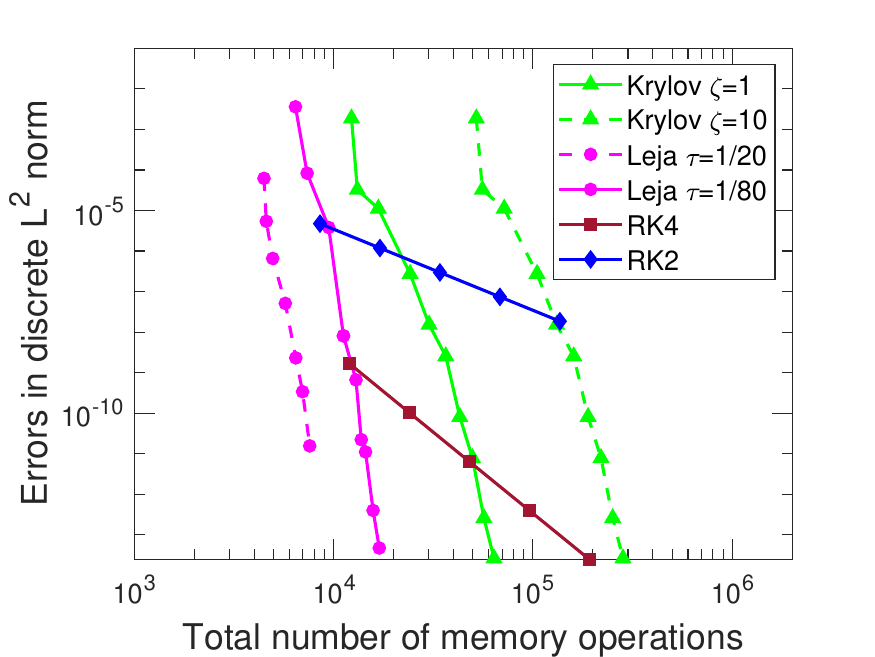}\label{fig4d}}					
\caption{The numerical results for the mixed problem, as defined in equation \eqref{problem1} and \eqref{eq12} are presented in the upper two panes. The figures illustrate the achieved accuracy and the required computational cost of the considered methods for evaluating $\exp(- (A_h+B_h))u_0$ as a function of the time step size. The tolerance is chosen such that the achieved accuracy for the exponential integrators is $10^{-4}$ (dashed-dotted line) and $10^{-7}$ (solid line), respectively. The pane~\ref{fig4cc} illustrates the function $\kappa$ in \eqref{eq12}. The pane~\ref{fig4d} depicts the achieved accuracy of the considered methods as a function of the computational cost for two different values of $\tau$ and $\zeta$. For Krylov methods the employed time step size is $\tau = \frac{1}{80}$.}\label{fig4}
\end{figure}

We proceed to investigate a mixed problem, in which diffusion dominates in one part of the domain and advection in the other. Setting
\begin{equation}\label{eq12}
\kappa(x) = \frac{33}{5120} + \frac{31}{5120} \tanh (20x-16)
\end{equation}
ensures that the desired properties are achieved as the value of $\kappa(x)$ spans from the advection-dominated case $\frac{1}{2560}$ to the diffusion-dominated case $\frac{1}{80}$ (see also Figure~\ref{fig4cc}). We choose $\tau = \frac{1}{10}$ as the maximum time step for the Krylov method and $\tau = \frac{1}{20}$ for the Leja method for evaluating the action of the matrix exponential.

Figure \ref{fig4} illustrates the numerical results. The approximation of $\exp(-(A_h + B_h))u_0$, in particular, using the Leja method is faster than explicit methods. This is particularly true for stringent tolerances.

\section{Numerical results for the compressible isothermal Navier--Stokes problem}\label{num2d}

Motivated by the promising experiments with the linear problem, this section is devoted to assessing the efficiency of exponential integrators in a more challenging nonlinear setting. To this end, we have selected the two-dimensional compressible isothermal Navier--Stokes problem for detailed investigation. As in the previous case, we will assess the performance of the selected exponential integrators in comparison to that of the standard explicit Runge--Kutta methods RK2 and RK4.

\subsection{Implementation}

\begin{figure}[b!]
\subfigure[$\rho(12,x,y)$]{\includegraphics[width=0.48\textwidth]{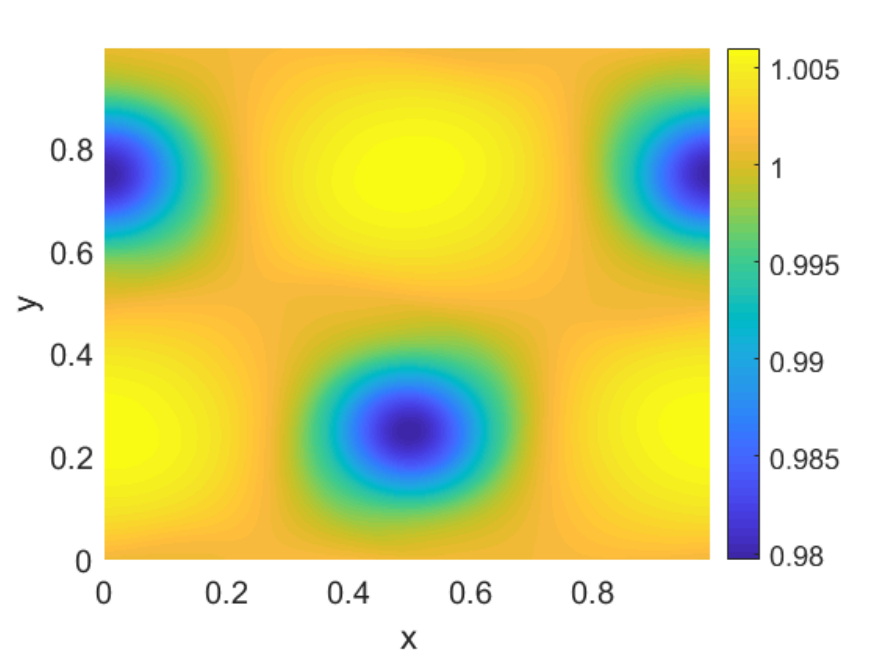}}
\subfigure[$u(12,x,y)$]{\includegraphics[width=0.48\textwidth]{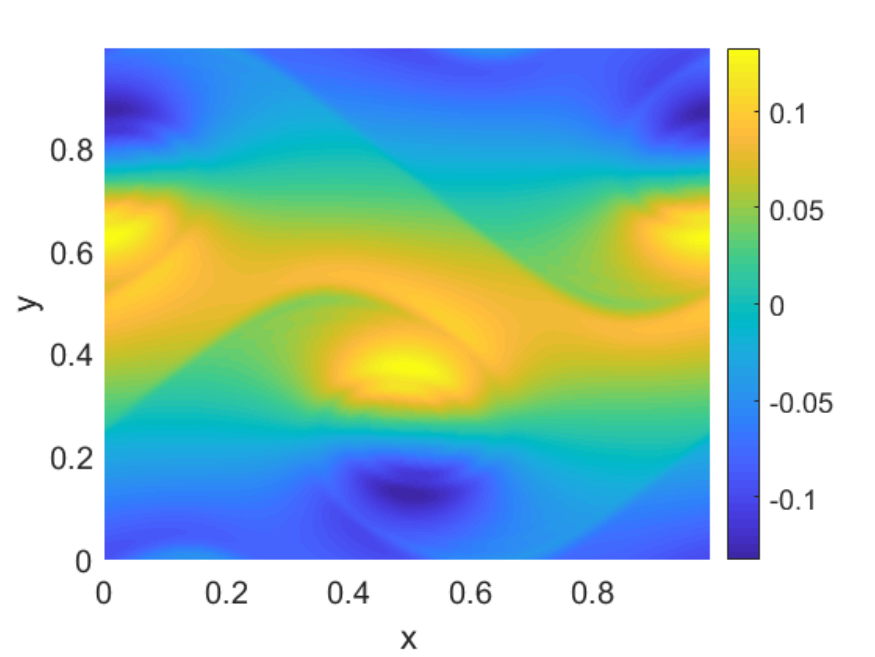}}
\subfigure[$v(12,x,y)$]{\includegraphics[width=0.48\textwidth]{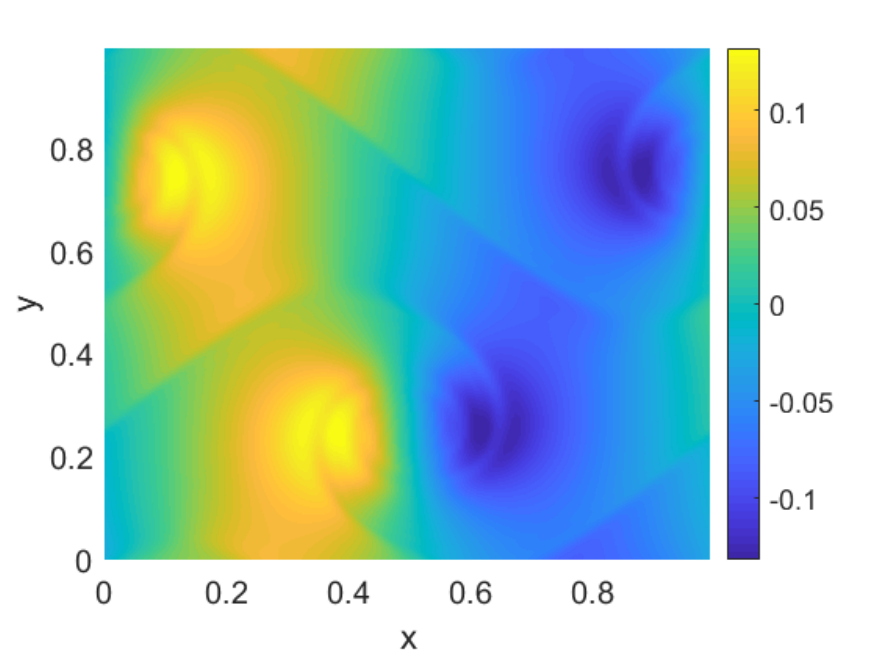}}
\subfigure[$\omega(12,x,y)$]{\includegraphics[width=0.48\textwidth]{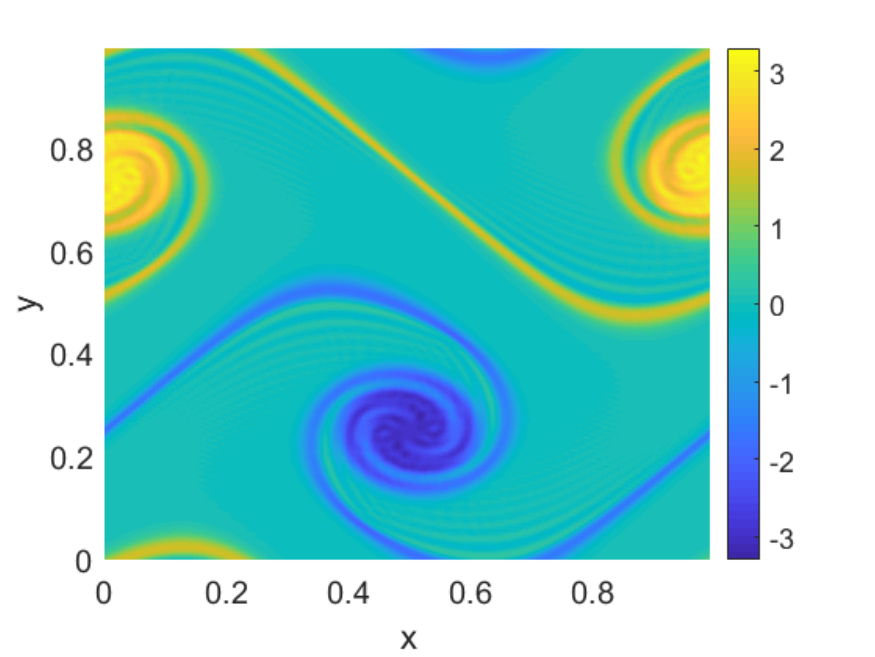}}	
\caption{The numerical simulation of the shear flow problem.}
\label{fig5}
\end{figure}

The discretization of \eqref{eq41} is performed using an equidistant mesh with $n = 160$ degrees of freedom in each direction. We take initial data corresponding to a shear flow, which is a classical example in computational fluid dynamics. A common setup (see, e.g., \cite{Lukas2022,EINKEMMER2014,LIU2000577}) is
\begin{align*}
		\rho(0, x, y) & =1,	\\
		u(0, x, y) &=\begin{cases}v_0 \tanh \bigl(\frac{y-\frac{1}{4}}{d}\bigr), & y \leq \frac{1}{2}, \\
			v_0 \tanh \bigl(\frac{\frac{3}{4}-y}{d}\bigr), &  y>\frac{1}{2},\end{cases} \\
		v(0, x, y) &=\delta \sin (2 \pi x),
\end{align*}
where $(x, y) \in[0,1]^{2}$, $v_{0}=0.1$, $d=1 / 30$, and $\delta=5 \cdot 10^{-3}$. From a physical perspective, the problem under consideration represents the motion of a fluid within the domain. Specifically, the fluid exhibits leftward movement at the top and bottom, while flowing in the opposite direction in the center of the domain. As time progresses, a small perturbation in the velocity field is amplified, leading to the formation of dynamic structures that are predominantly characterized by vortices. The (numerical) solution, obtained at the final time $t = 12$, is presented in Fig.~\ref{fig5}.

It is important to note that, in contrast to the linear problem, the Jacobian matrix is now changing at each time step. This is not an issue for explicit methods, which require only the evaluation of the right-hand side of the problem. For exponential integrators, however, the issue of a changing Jacobian requires particular care and attention. Krylov methods must build a new sequence of subspaces and interpolation methods, such as Leja interpolation, require some knowledge of the spectrum of the Jacobian. For further details on how to estimate the spectrum for cases where the coefficients of the matrix are explicitly known (the case considered in this paper), we refer the reader to the reference \cite{CaKOR16}. For matrix-free implementations a method based on the power iteration has been discussed in \cite{pub.1146932560}.

\subsection{Numerical investigations}

\begin{figure}[h!]
\subfigure{\includegraphics[width=0.48\textwidth]{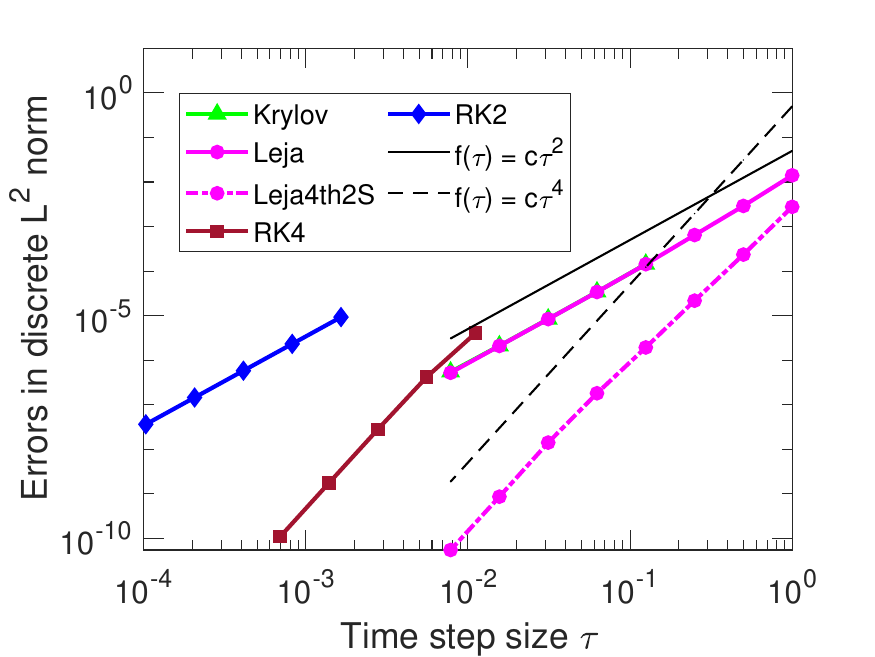}}
\subfigure{\includegraphics[width=0.48\textwidth]{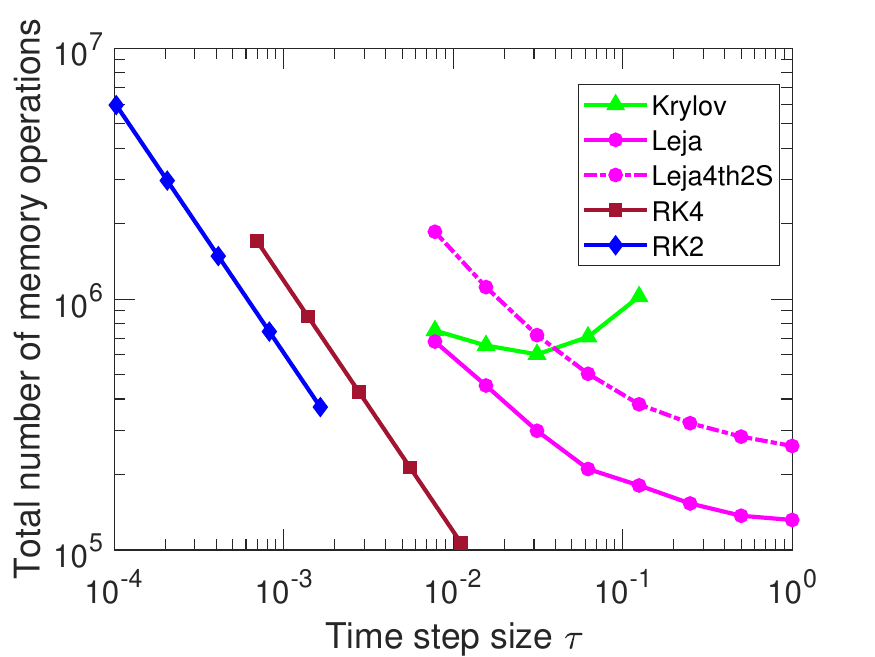}}	
\begin{center}
\subfigure{\includegraphics[width=0.48\textwidth]{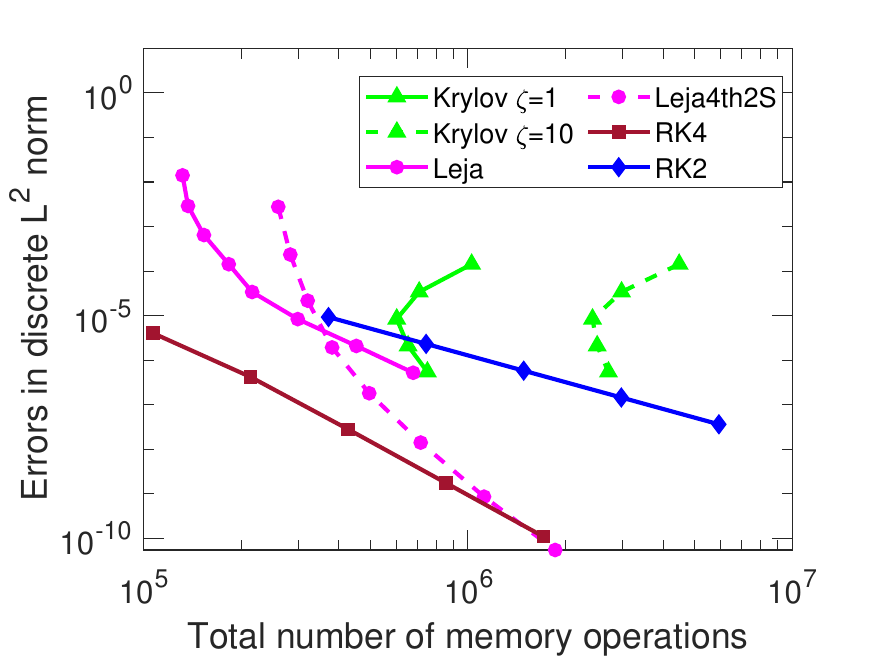}}	
\end{center}
\caption{The numerical results for the two-dimensional compressible isothermal Navier--Stokes problem \eqref{eq41} with Reynolds number $10^5$ at time $t=12$. The two figures at the top illustrate the global error and the computational cost as functions of the time step size for the various methods under consideration. For a given time step size, the errors for the exponential Rosenbrock--Euler method, when using the Krylov method or Leja interpolation, are similar. The bottom figure shows a comparison of the cost of the various methods, wherein the exponential methods use the time step sizes $\tau = 2^{-m}$, $0\le m\le 7$.}\label{fig3}
\end{figure}

All simulations are conducted up to the final time $t = 12$. For a given Reynolds number Re, the value of $\nu$ is determined by the formula $\nu = v_0 / \text{Re}$. We choose $\text{Re} = 10^5$, which yields a kinematic viscosity of $\nu = 10^{-6}$. As a consequence, the problem is advection dominated. We further employ $\tau = 1$ as the maximum time step size for the exponential integrators. Figure \ref{fig3} presents the work-precision diagram for the employed methods. Note that Leja4th2S denotes the scheme $\mathtt{exprb42}$ using the Leja method. For computational efficiency, the action of the matrix function $\varphi_1$ is computed directly in the first stage of $\mathtt{exprb42}$ and in the exponential Rosenbrock--Euler method. For updating the solution of $\mathtt{exprb42}$, however, we use the augmented matrix as described in \cite{doi:10.1137/100788860}. The reference solution is computed with RK4 using a sufficiently small step size.
	
Exponential integrators either based on the Krylov or the Leja method yield similar advantages to those observed in the linear advection-dominated case. The methods allow for significantly larger step sizes than those employed by RK2 and RK4, up to a factor of 600 and 90 times, respectively. Furthermore, the exponential integrators based on the Leja method exhibit superior computational efficiency compared to the RK2 method and are comparable to the RK4 method, as shown by the total number of required memory operations.

In regard to second-order integrators (specifically, the comparison of the exponential Rosenbrock--Euler method with RK2), the Leja-based method can achieve cost reductions of approximately 1.3 while maintaining a similar level of accuracy within the range of $10^{-5}$. In contrast, the fourth-order scheme $\mathtt{exprb42}$ is approximately three and a half times more expensive than the RK4 method, while achieving an accuracy of $10^{-5}$. However, if an accuracy of $10^{-3}$ is prescribed, which is sufficient for most practical applications, the Leja-based method is only about twice as expensive as RK4. Nevertheless, the $\mathtt{exprb42}$ scheme is competitive with RK4 when using smaller time steps (i.e.~for stringent tolerances), since the cost is not doubled when the step size is halved, in contrast to RK4.

\subsection{Comparison of the various methods}

The exponential Rosenbrock--Euler scheme has been demonstrated to achieve an accuracy level of approximately $10^{-3}$ even for large time steps, which is sufficient for the majority of applications. Furthermore, this level of accuracy requires much less computational effort than that required by RK2. Note that in this regime both RK2 and RK4 are constrained by stability conditions. While Krylov-based methods are expensive for large time step sizes, they are more effective than RK2 when small time steps are employed. As can be seen, the exponential Rosenbrock--Euler method (based on the Leja method) consistently outperforms the RK2 method across all tolerance regimes, from medium to extremely stringent tolerances (this is most likely helped by the bad linear stability properties of RK2 on the imaginary axis).

As with the linear problem, variations in the parameter $\zeta$ primarily affect Krylov methods, while the other methods are less affected. This makes the Krylov methods expensive when implemented on a supercomputer. The Leja-based exponential fourth-order integrator lags behind the RK4 method for large time steps. Nevertheless, it becomes competitive at stringent tolerances.

\section{Conclusion}\label{sectionconcludesionpaper1}

We have studied the efficiency of exponential integrators for advection-dominated problems. While Leja-based exponential integrators perform similar to explicit schemes, in contrast to diffusion-dominated problems, no significant improvement in performance can be observed. However, in situations where parts of the domain are diffusion dominated, exponential integrators are able to significantly outperform explicit schemes. While exponential integrators are able to take large time steps, even in the advection-dominated case, also a relative large number of iterations is required in each time step. An advantage of exponential integrators in the linear case is that they offer high precision at little additional computational cost.

In all the problems considered we observe that Leja-based methods outperform Krylov iterations. In particular, Leja methods are attractive on modern supercomputers as they do not require the computation of inner products.

\bibliography{references}

\appendix

\section{Appendix}\label{Appendix1}

This appendix contains some technical details of Section \ref{discussNavierstokes}. Recall that we have
\begin{equation*}
	U = \left[\begin{array}{l}
		\rho \\
		u \\
		v
	\end{array}\right], \quad F(U) = \left[\begin{array}{c}
		F_1(U) \\
		F_2(U) \\
		F_3(U)
	\end{array}\right] =
	\left[\begin{array}{c}
		-\partial_x (\rho u) - \partial_y (\rho v) \\
		- u \partial_{x} u - v \partial_{y} u - \frac{1}{\rho} \partial_{x} \rho + \nu (\partial_{xx}u + \partial_{yy}u) \\
		-u \partial_{x} v-v \partial_{y} v - \frac{1}{\rho} \partial_{y} \rho +\nu\left(\partial_{x x} v+\partial_{y y} v\right)
	\end{array}\right].
\end{equation*}
By differentiating the functions $F_1$, $F_2$, and $F_3$ with respect to the variables $\rho$, $u$, and $v$ respectively, we obtain the (actions of the) components of the Jacobian
\begin{gather*}
\frac{\partial F_{1}}{\partial \rho} (U)(w) = - \partial_x (uw) - \partial_y (vw),\quad \frac{\partial F_{1}}{\partial u} (U)(w) = - \partial_x (\rho w), \quad \frac{\partial F_{1}}{\partial v} (U)(w) = - \partial_y (\rho w),\\
\frac{\partial F_{2}}{\partial \rho} (U)(w) = w\frac{1}{\rho^2} \partial_x \rho -\frac{1}{\rho} \partial_x w, \quad \frac{\partial F_{2}}{\partial u} (U)(w) = - \partial_x( u w) - v \partial_y w + \nu (\partial_{xx} w + \partial_{yy} w),\\
\frac{\partial F_{2}}{\partial v} (U)(w) =  - w \partial_y u, \quad \frac{\partial F_{3}}{\partial \rho} (U)(w) = w\frac{1}{\rho^2} \partial_y \rho -\frac{1}{\rho} \partial_y w ,\\
\frac{\partial F_{3}}{\partial u} (U)(w) =  - w \partial_x v, \quad  \frac{\partial F_{3}}{\partial v} (U)(w) = - u \partial_x w - \partial_y( v w )  + \nu (\partial_{xx} w + \partial_{yy} w) .
\end{gather*}
For the spatial discretization of $U$, we use again standard second-order finite differences on an equidistant mesh on $\Omega = [0,1]$ with mesh width $h$. Let $n$ denote the number of grid points on each axis, and let the matrices $A_h$ and $B_h$ denote the finite difference discretizations of $\nu \partial_{xx}$ (or $\nu \partial_{yy}$) and $\partial_{x}$ (or $\partial_{y}$) for periodic boundary conditions. The matrix formulation of the action of the Jacobian matrix applied to the vector $[w_1,w_2,w_3]^T$ is as follows
\begin{align*}
	J_{11}(U)w_1 &+ J_{12}(U)w_2 + J_{13}(U)w_3 = \\
	&  -({I}  \otimes {B_h}) (u \odot w_1) - ({B_h} \otimes {I}) (v \odot w_1) -({I} \otimes {B_h}) (\rho \odot w_2) -({B_h} \otimes {I}) (\rho \odot w_3), \\[2pt]
	J_{21}(U)w_1 &+  J_{22}(U)w_2 + J_{23}(U)w_3 = \\
	& w_1 \odot \tfrac{1}{\rho^2} \odot (({I} \otimes {B_h})\rho) - \tfrac{1}{\rho} \odot (({I} \otimes {B_h})w_1) -({I} \otimes {B_h})(u \odot w_2) - v\odot (({B_h} \otimes {I})w_2) \\
	&+  (  ({I} \otimes {A_h}) + ({A_h} \otimes {I})  )w_2- w_2 \odot (({B_h} \otimes {I})u ),\\[2pt]
	J_{31}(U)w_1 &+ J_{32}(U)w_2 + J_{33}(U)w_3 = \\
	& w_1 \odot \tfrac{1}{\rho^2} \odot (({B_h} \otimes {I})\rho) -\tfrac{1}{\rho}\odot (({B_h} \otimes {I})w_1)  - w_3 \odot (({I} \otimes {B_h})v) \\
	& - u\odot (({I} \otimes {B_h})w_3) -({B_h} \otimes {I}) (v \odot w_3)  +  (({I} \otimes {A_h}) + ({A_h} \otimes {I})   )w_3.
\end{align*}
Here,
$$
u \odot v = [u(1)v(1), u(2)v(2), \cdots, u(N)v(N)]^T
$$
denotes the component-wise product of two vectors and $N=n^2$ represents the length of the occurring vectors. Explicit Runge--Kutta methods only evaluate the vector $F(U)$. Each component of this vector is computed as follows
\begin{gather*}
F_1(U) = - ({I} \otimes {B_h})(\rho \odot u) - ({B_h} \otimes {I}) (\rho \odot v),\\
F_2(U) = - u \odot ({I} \otimes {B_h}) u - v \odot ({B_h} \otimes {I})u - \tfrac{1}{\rho} \odot (({I} \otimes {B_h} ) \rho) +  ({I} \otimes {A_h} +  {A_h} \otimes {I}) u, \\
F_3(U) = - u\odot ({I} \otimes {B_h}) v - v\odot ({B_h} \otimes {I})v - \tfrac{1}{\rho} \odot (({B_h} \otimes {I}) \rho) +  ({I} \otimes {A_h} + {A_h} \otimes {I}) v .
\end{gather*}

\end{document}